\documentclass [12pt]{article}
\usepackage {multicol}

\usepackage{amssymb,amsmath}
\usepackage{mathrsfs}
\usepackage{hyperref}
\usepackage{graphicx}
\usepackage{color}
\usepackage{caption}

\begin{document}

\title{\Large\bf  Approximative Theorem of Incomplete Riemann-Stieltjes Sum of Stochastic Integral}

\author{\normalsize  Jingwei Liu\footnote{\scriptsize{Corresponding Author: Email: liujingwei03@tsinghua.org.cn (J.W Liu)}} \\
{\scriptsize{( School of Mathematics and System Sciences, Beihang University, Beijing,100083,P.R.China)}}
}

\date{\normalsize March 14,2018}
\maketitle

\textbf{Abstract:} The approximative theorems of incomplete Riemann-Stieltjes sums of Ito stochastic integral, mean square integral and Stratonovich stochastic integral with respect to Brownian motion are investigated. Some sufficient conditions of incomplete Riemann-Stieltjes sums approaching stochastic integral are developed, which establish the alternative ways to converge stochastic integral. And, Two simulation examples of incomplete Riemann-Stieltjes sums about Ito stochastic integral and Stratonovich stochastic integral are given for demonstration.

\textbf{Keywords:} Ito stochastic integral; Stratonovich stochastic integral; Mean square integral; Definite integral; Brownian motion ;Riemann-Stieltjes sum; Convergence in mean square; Convergence in probability

\section{Introduction}

Ito stochastic integral and  Stratonovich stochastic integral take important roles in stochastic calculus and stochastic differential equation. Both of them can be expressed in the limit of Riemann-Stieltjes sum [1-4], though they can be modeled by martingale and semimartingale theory[5-13]. [14] proposes the deleting item  Riemann-Stieltjes sum theorem of definite integral, which can also be called incomplete Riemann-Stieltjes sum. As definitions of martingale and semimartingale are still relative to integral including definite integral, and the stochastic Riemann-Stieltjes sum limit provides convenient way for numeral solution, the most important, stochastic integral and stochastic differential equation are mutual determined , is there alternative way of stochastic integral to describe the stochastic integral and stochastic differential equation convergence in mean square limit or in probability sense? In fact, [14] has revealed that definite integral is not uniquely determined by differential equation according to deleting item Riemann-Stieltjes sum theorem. Furthermore, there is another stochastic integral, mean square integral [15-18], is involved in stochastic calculus. Motivated by above reasons, we concentrate on the limits of incomplete Riemann--Stieltjes sums of Ito stochastic integral, Stratonovich stochastic integral and mean square integral.

The rest of the paper is structured as follows: Section 2 briefly reviews definitions of Ito stochastic integral, mean square integral and Stratonovich stochastic integral. Section 3 constructs incomplete Riemann-Stieltjes sum approximating theorems of Ito stochastic integral, mean square integral and Stratonovich stochastic integral. Section 4 extends incomplete Riemann-Stieltjes sum approximating theorem of stochastic integral. Simulations of incomplete Riemann-Stieltjes sum approximating are demonstrated in Section 5 and discussion is given in Section 6.

\section{Brief review of stochastic integral}

Let $(\Omega,\mathcal{F},(\mathcal{F}_t)_{t\in [0,T]},P)$ be a probability space where $(\mathcal{F}_t)_{t\in [0,T]}$ is a complete right continuous filtration generated by 1-dimensional standard Brownian motion (or Wiener process) $(B_t)_{t\in [0.T]}$ and each $\mathcal{F}_t$ containing all $P-$null sets of $\mathcal{F}_0$. A stochastic process $\Phi(\omega,t)$ is a predictable or non-anticipating process if $\Phi(\omega,t) \in (\mathcal{F}_t)$ for each $t$, where $\omega \in \Omega$, which is also called adapted stochastic process. And, denote a stochastic process $X(\cdot,t)$ as $X(t)=X_t \triangleq X(\omega,t)$, $\Phi_t\triangleq \Phi(\omega,t)$. Denote $\mathop{l.i.m}$ as limit in mean square.

\textbf{Definition 1} Suppose that $\Phi(\cdot,t)$ is a non-anticipating stochastic process with respect to the Brownian motion $B_{t}$. Let  $\bigtriangleup_{n}
= \{t_{0},t_{1},\cdots,t_{n}\}$ be the partition : $0=t_{0}<t_{1}<\cdots<t_{n}=T$ of $[0,T]$, $\vartriangle t_i= (t_{i+1}- t_{i})$, $\|\bigtriangleup_{n}\|= \displaystyle \max_{0\leq j\leq n-1} \vartriangle t_i $. The It$\hat{o}$ stochastic integral of the stochastic process $\Phi(\cdot,t)$ is denoted by
\begin{equation}
\displaystyle \int_{0}^{T} \Phi_{t} dW_{t}=\displaystyle \mathop{l.i.m}\limits_{\|\bigtriangleup_{n}\| \longrightarrow 0} \sum_{i=0}^{n-1} \Phi(B_{t_{i}},t_{i}) (B_{t_{i+1}}-B_{t_{i}})
\end{equation}

\textbf{Definition 2} (Mean-square integration) Let $X(\cdot,t)$ be a second-order stochastic process on $\Omega\times [0,T]$, which is $E[X^2(\cdot,t)]<+\infty$. For any partition :$0=t_{0}<t_{1}<\cdots <t_{n}=T $, $\|\bigtriangleup_{n}\|=\displaystyle \max_{0\leq i\leq n-1} (t_{i+1}- t_{i})$, and, for any $u_i\in[t_{i},t_{i+1}]$. The mean--square integral $\int_{0}^{T} X(t)dt$ is defined as the limit in mean square
\begin{equation}
\displaystyle \int_{0}^{T} X(t) dt=\displaystyle \mathop{l.i.m}\limits_{\|\bigtriangleup_{n}\|\longrightarrow 0} \sum_{i=0}^{n-1} X(u_{i})(t_{i+1}- t_{i})
\end{equation}
provided that the limit exists for arbitrary partition and selected points $u_i$. $X(t)$ is called mean-square integrable on $[0,T]$.

If  $\displaystyle \Phi(B_{t_{i}},t_{i})$ is replaced by $\displaystyle \Phi(\frac{B_{t_{i}}+B_{t_{i+1}}}{2},t_{i})$ in It$\hat{o}$ stochastic integral, it will be Stratonovich stochastic integral.

\textbf{Definition 3} Suppose that $\Phi(\cdot,t)$ is a non-anticipating stochastic process with respect to the Brownian motion $B_{t}$. Let  $\bigtriangleup_{n}
= \{t_{0},t_{1},\cdots,t_{n}\}$ be the partition : $0=t_{0}<t_{1}<\cdots<t_{n}=T$ of $[0,T]$, $\vartriangle t_i= (t_{i+1}- t_{i})$, $\|\bigtriangleup_{n}\|=  \max \vartriangle t_i $. The Stratonovich integral of the stochastic process $\Phi(\cdot,t)$ is denoted by
\begin{equation}
\displaystyle \int_{0}^{T} \Phi_{t}\circ dW_{t}=\displaystyle \mathop{l.i.m}\limits_{\|\bigtriangleup_{n}\| \longrightarrow 0} \sum_{i=0}^{n-1} \Phi(\frac{B_{t_{i}}+B_{t_{i+1}}}{2},t_{i}) (B_{t_{i+1}}-B_{t_{i}})
\end{equation}

In fact, Stratonovich stochastic integral has an alternative formula as follows [3,4], it can be achieved by derivation while $\Phi(x,t)$ is 1-order differential with $x$ [8].
\begin{equation}
\displaystyle \int_{0}^{T} \Phi_{t}\circ dW_{t}=\displaystyle \mathop{l.i.m}\limits_{\|\bigtriangleup_{n}\| \longrightarrow 0} \sum_{i=0}^{n-1}  \frac{\Phi_{t_{i+1}}+\Phi_{t_{i}}}{2} (B_{t_{i+1}}-B_{t_{i}})
\end{equation}

If function $\Phi(x,t)$ continuous in $t$, having the continuous derivative $\displaystyle \frac{\partial \Phi(x,t)}{\partial x} $. There is a relationship among It$\hat{o}$ stochastic integral, mean square integral and Stratonovich stochastic integral as follows.
\begin{equation}
\displaystyle \int_{0}^{T} \Phi(B(t),t)\circ  dB(t)=\int_{0}^{T} \Phi(B(t),t) dB(t) + \displaystyle \frac{1}{2} \int_{0}^{T}\frac{\partial \Phi(B(t),t)}{\partial B(t)}dt
\end{equation}


Denote $J = \{0,1,2, \cdots,n-1 \}$. For a fixed natural number $K \in N^{+} (0 < K < n)$, denote $J_K =
\{i_1 , \cdots ,i_K \} \subset J$, and $J\setminus J_{K}=\{0,1,2,\cdots,n-1 \}\setminus \{i_1 , \cdots ,i_K \} $.

We will investigate the limit of incomplete Riemann-Stieltjes sums in mean square convergence,
\begin{equation}
\displaystyle  \sum_{J\setminus J_{K} } \Phi(B_{t_{i}},t_{i}) (B_{t_{i+1}}-B_{t_{i}})
\end{equation}
\begin{equation}
\displaystyle  \sum_{J\setminus J_{K} } X(u_{i}) (t_{i+1}-t_{i})
\end{equation}
and,
\begin{equation}
\displaystyle \sum_{J\setminus J_{K}} \Phi(\frac{B(t_{j})+B(t_{j+1})}{2},t_{j}) [B(t_{j+1})-B(t_{j})]
\end{equation}

In [1], Ito Stochastic integral is originally defined on $B(\omega,t)$ with no moving discontinuity, such that $\mathbb{E}[B(\omega,t)-B(\omega,s)]=0$ , $\mathbb{E} [B(\omega,t)-B(\omega,s)]^2=|t-s|$. Our investigation on stochastic integral with incomplete Riemann-Stieltjes sum will show the asymptotical limit behavior of discontinuity subintervals on continuous path of Brownian motion, which is, for any partition of $[s,t]$, $\displaystyle  \sum_{J\setminus J_{K} } \mathbb{E} [B_{t_{i+1}}-B_{t_{i}}]=0$, and $\displaystyle  \sum_{J\setminus J_{K} } \mathbb{E} [B_{t_{i+1}}-B_{t_{i}}]^2 \longrightarrow |t-s|$,  as $n\longrightarrow +\infty$ . Our motivation is also different from the $\delta-$fine belated partial division of $[s,t]$ in [19]. Here, $[s,t]=[0,T]$.

Two convergence properties, convergence in mean square ($\mathop{l.i.m}$ or $\stackrel{m.s}{\Longrightarrow }$ ) and convergence in probability ($\stackrel{P}{\longrightarrow}$), and two kinds of inequalities are involved in the discussion [20,21].

\textbf{Lemma 1. } Suppose that $X_n\stackrel{m.s}{\Longrightarrow } X$, $Y_n\stackrel{m.s}{\Longrightarrow } Y$. For $\forall a, b \in \mathbf{R}$, then
\begin{equation}
 aX_n+ bY_n \stackrel{m.s}{\Longrightarrow } aX+ bY.
\end{equation}

\textbf{Lemma 2. }Suppose that $X_n\stackrel{P}{\longrightarrow } X$, $Y_n\stackrel{P}{\longrightarrow } Y$. For $\forall a, b \in \mathbf{R}$, then
\begin{equation}
 aX_n+ bY_n \stackrel{P}{\longrightarrow } aX+ bY.
\end{equation}

\textbf{Lemma 3. } (H$\ddot{o}$lder Inequality) Suppose that $\xi,\eta$ are two random variables on $(\Omega,\mathcal{F},P)$, for any real numbers $p,q \in R$ satisfying $1<p,q<+\infty$, and $\displaystyle \frac{1}{p}+\frac{1}{q}=1$,
\begin{equation}
 E|\xi\eta|\leq (E^{1/p}|\xi|^p)(E^{1/q}|\eta|^q).
\end{equation}
If $E|\xi|^p<+\infty$ and $E|\eta|^q<+\infty$, then ``='' holds only when $\xi=0$ a.s. or $\eta=0$ a.s. ,or there exists constant $C$ so that $|\xi|^p=C |\eta|^q $ a.s. .

\textbf{Lemma 4. } (Minkowski Inequality) Suppose that $\xi,\eta$ are two random variables on $(\Omega,\mathcal{F},P)$,

1) For any real number $p\geq 1$ ,
\begin{equation}
 E^{1/p}|\xi+\eta|^p \leq E^{1/p}|\xi|^p + E^{1/p}|\eta|^p.
\end{equation}
If $E(|\xi|^p +|\eta|^p)<+\infty$, then ``='' holds only when $p>1$, $\xi=0$ a.s. or $\eta=0$ a.s. ,or there exists constant $C>0$ so that $\xi=C\eta$ a.s.;when $p=1$, $\xi\eta\geq 0$ a.s. .

2) For $0<p<1$,
\begin{equation}
 E|\xi+\eta|^p \leq E|\xi|^p + E|\eta|^p).
\end{equation}
If $E(|\xi|^p +|\eta|^p)<+\infty$, then ``=''  holds if and only if $\xi\eta=0$ a.s. .

\section{Incomplete Riemann-Stieltjes Sum approximating of Stochastic integral}

Since Brownian motion is of unbounded variation, the convergence of Riemann-Stieltjes sum of Ito stochastic integral is not in usual Riemann-Stieltjes  sense, but something like Lebesgue sense. It is first defined on a class of simple functions, then extended to a large class of functions with isometry [1,6]. The following notation is from [13].

Denote
\begin{equation}
\begin{array}{ll}
\mathscr{L}_{0}&=\{\Phi: \Phi=f_{0}(\omega) I_{0}(t)+\displaystyle \sum_{k=1}^{\infty} f_{k}(\omega) I_{(t_{k},t_{k+1}]}(t),
                0=t_{0}<t_{1}<\cdots<t_{n} \uparrow \infty . \\
               & \sup_{\substack{0\leq k <\infty}} \{ f_{k} \} < \infty , f_{0} \in \mathscr{F}_{0},  f_{k} \in \mathscr{F}_{k} \}
\end{array}
\end{equation}

\begin{equation}
\begin{array}{ll}
\mathscr{L}_{p}&=\{\Phi: \Phi=(\Phi_{t})_{t\geq 0} , \Phi_{t}(\omega) \equiv \Phi(\omega,t)  \mbox{ is measurable and }
(\mathscr{F}_{t}) \mbox{ adapted} \\
& \mbox{process}, \mbox{and for } \forall   T , \mathbb{E} \int_{0}^{T} |\Phi_{t}|^{p} dt <\infty  \} ,  p\geq 1.
\end{array}
\end{equation}

\begin{equation}
\begin{array}{ll}
\mathscr{L}_{p,T}&=\{\Phi: \Phi=(\Phi_{t})_{t\geq 0} , \Phi_{t}(\omega) \equiv \Phi(\omega,t)  \mbox{ is measurable and }
(\mathscr{F}_{t}) \mbox{ adapted} \\
& \mbox{process}, \mathbb{E} \int_{0}^{T} |\Phi_{t}|^{p} dt <\infty  \} ,  p\geq 1.
\end{array}
\end{equation}

\begin{equation}
\begin{array}{ll}
L_{p}&=\{\Phi: \Phi(\omega,t), \omega\in\Omega, t\in[0,T], \mathbb{E} \int_{0}^{T} |\Phi_{t}|^{p} dt <\infty  \} ,  p\geq 1.
\end{array}
\end{equation}

\begin{equation}
\begin{array}{ll}
\mathscr{L}_{p}^{loc}&=\{\Phi: \Phi=(\Phi_{t})_{t\geq 0} , \Phi_{t}(\omega) \equiv \Phi(\omega,t)  \mbox{ is measurable and }
(\mathscr{F}_{t}) \mbox{ adapted} \\
& \mbox{process}, \mbox{and for } \forall   T , \int_{0}^{T} |\Phi_{t}|^{p} dt <\infty  \} ,  p\geq 1.
\end{array}
\end{equation}

\begin{equation}
\begin{array}{ll}
\mathscr{L}_{p,T}^{loc}&=\{\Phi: \Phi=(\Phi_{t})_{t\geq 0} , \Phi_{t}(\omega) \equiv \Phi(\omega,t)  \mbox{ is measurable and }
(\mathscr{F}_{t}) \mbox{ adapted} \\
& \mbox{process}, \int_{0}^{T} |\Phi_{t}|^{p} dt <\infty  \} ,  p\geq 1.
\end{array}
\end{equation}

The classical definition of Ito stochastic integral is defined on $\mathscr{L}_{0}$,$\mathscr{L}_{2}$ and $\mathscr{L}_{2}^{loc}$ [13].

\subsection{Incomplete Riemann-Stieltjes Sum approximating of Ito Stochastic integral}

\textbf{Lemma 5. } Suppose that $\Phi(\cdot,t)$ is  nonanticipative and  $\Phi(\cdot,t)\in \mathscr{L}_{0}$, then
\begin{equation}
 \mathbb{E}[ \sum_{ \substack{j\in J_{K}} }  \Phi(B_{t_{j}},t_{j}) (B_{t_{j+1}}-B_{t_{j}})]= 0 .
\end{equation}
\begin{equation}
 \mathbb{E} | \sum_{ \substack{j\in J_{K}} }  \Phi(B_{t_{j}},t_{j}) (B_{t_{j+1}}-B_{t_{j}}) -0 |^{2}\longrightarrow 0,  \  \  as \ \ \|\bigtriangleup_{n}\| \longrightarrow 0.
\end{equation}
And,
\begin{equation}
  \sum_{ \substack{j\in J_{K}} }  \Phi(B_{t_{j}},t_{j}) (B_{t_{j+1}}-B_{t_{j}}) \stackrel{P}{\longrightarrow 0},\ \ as \ \ \|\bigtriangleup_{n}\| \longrightarrow 0.
\end{equation}

\textbf{Proof. } (1) Since $ \Phi_{t}$ is nonanticipative, and $B_{t}$ is incremental independent process.
\begin{equation}
\displaystyle \mathbb{E} \sum_{ \substack{j\in J_{K}} }  \Phi(B_{t_{j}},t_{j}) (B_{t_{j+1}}-B_{t_{j}})
=\displaystyle  \sum_{ \substack{j\in J_{K}} } \mathbb{E} \Phi(B_{t_{j}},t_{j}) \mathbb{E}(B_{t_{j+1}}-B_{t_{j}})=0
\end{equation}

(2) Using above conclusion, we obtain
\begin{equation}
\begin{array}{ll}
&\displaystyle \mathbb{E} | \sum_{ \substack{j\in J_{K}} }  \Phi(B_{t_{j}},t_{j}) (B_{t_{j+1}}-B_{t_{j}}) |^2\\
   =&\displaystyle \mathbb{E}  \sum_{ \substack{j\in J_{K}} }  |\Phi(B_{t_{j}},t_{j}) (B_{t_{j+1}}-B_{t_{j}}) |^2  \\
   =&\displaystyle  \sum_{ \substack{j\in J_{K}} } \mathbb{E} |\Phi(B_{t_{j}},t_{j}) (B_{t_{j+1}}-B_{t_{j}}) |^2\\
   =&\displaystyle  \sum_{ \substack{j\in J_{K}} } \mathbb{E} |\Phi^2(B_{t_{j}},t_{j}) (B_{t_{j+1}}-B_{t_{j}})^2 | \\
\leq & \displaystyle  \sum_{ \substack{j\in J_{K}} } \mathbb{E} \sup_{\substack{j\in J_{K}}} \{\Phi^2(B_{t_{j}},t_{j})\} (B_{t_{j+1}}-B_{t_{j}})^2\\
\leq & \displaystyle  \sum_{ \substack{j\in J_{K}} }  \sup_{\substack{j\in J_{K}}} \{\Phi^2(B_{t_{j}},t_{j})\}\mathbb{E} (B_{t_{j+1}}-B_{t_{j}})^2  \\
\leq &\displaystyle  \sum_{ \substack{j\in J_{K}} }  \sup_{\substack{j\in J_{K}}} \{\Phi^2(B_{t_{j}},t_{j})\} {\vartriangle t_{j}}\\
\leq &\displaystyle  \sup_{\substack{j\in J_{K}}} \{\Phi^2(B_{t_{j}},t_{j})\} \sum_{ \substack{j\in J_{K}}} \{\vartriangle t_{j}\} \\
\leq &\displaystyle  \sup_{\substack{j\in J_{K}}} \{\Phi^2(B_{t_{j}},t_{j})\} K \max{ \{\vartriangle t_{j}\}} \longrightarrow 0, \  \ as \ \ {\max{ \{\vartriangle t_{i}\}}\longrightarrow 0} \\
\end{array}
\end{equation}
Hence,
\begin{equation}
\displaystyle \mathop{l.i.m}\limits_{\|\bigtriangleup_{n}\| \longrightarrow 0 }  \sum_{ \substack{j\in J_{K}} }  \Phi(B_{t_{j}},t_{j}) (B_{t_{j+1}}-B_{t_{j}}) =0
\end{equation}
Then,
\begin{equation}
\displaystyle  \sum_{ \substack{j\in J_{K}} }  \Phi(B_{t_{j}},t_{j}) (B_{t_{j+1}}-B_{t_{j}})\stackrel{P}{\longrightarrow 0},\ \ as \|\bigtriangleup_{n}\| \longrightarrow 0.
\end{equation}

\hfill $\Box$

\textbf{Theorem 1. } Suppose that  $\Phi(\cdot,t)\in \mathscr{L}_{0} \bigcup \mathscr{L}_{2} \bigcup \mathscr{L}_{2}^{loc}$, then
\begin{equation}
\displaystyle \mathop{l.i.m}\limits_{\|\bigtriangleup_{n}\| \longrightarrow 0 }  \sum_{ \substack{j\in J\setminus J_{K}} }  \Phi(B_{t_{j}},t_{j}) (B_{t_{j+1}}-B_{t_{j}})= \int_0^T \Phi_{t} d B_{t}.
\end{equation}

\textbf{Proof. } As in [13], the classical definition of Ito stochastic process is defined on $\mathscr{L}_{0}$,$\mathscr{L}_{2}$ and $\mathscr{L}_{2}^{loc}$. We prove the conclusion step by step in a concise way, the detail proof will refer to [13].

(1) For $ \Phi_{t} \in \mathscr{L}_{0}$, we first discuss $\displaystyle \sum_{ \substack{j\in J_{K}} } \Phi(B_{t_{j}},t_{j}) (B_{t_{j+1}}-B_{t_{j}}) $ . \\

According to Lemma 5 , we obtain

\begin{equation}
\displaystyle \mathop{l.i.m}\limits_{\|\bigtriangleup_{n}\| \longrightarrow 0 }  \sum_{ \substack{j\in J_{K}} }  \Phi(B_{t_{j}},t_{j}) (B_{t_{j+1}}-B_{t_{j}}) =0
\end{equation}

According to Lemma 1, we obtain
\begin{equation}
\begin{array}{ll}
&\displaystyle \mathop{l.i.m}\limits_{\|\bigtriangleup_{n}\| \longrightarrow 0}  \sum_{ \substack{j\in J\setminus J_{K}} }  \Phi(B_{t_{j}},t_{j}) (B_{t_{j+1}}-B_{t_{j}}) \\
&= \displaystyle \displaystyle \mathop{l.i.m}\limits_{\|\bigtriangleup_{n}\| \longrightarrow 0}  \sum_{ \substack{j\in J} } \Phi(B_{t_{j}},t_{j}) (B_{t_{j+1}}-B_{t_{j}}) -  \mathop{l.i.m}\limits_{\|\bigtriangleup_{n}\| \longrightarrow 0} \sum_{ \substack{j\in J_{K}} }  \Phi(B_{t_{j}},t_{j}) (B_{t_{j+1}}-B_{t_{j}})\\
&= \int_0^T \Phi_{t} d B_{t} -0 = \int_0^T \Phi_{t} d B_{t}.
\end{array}
\end{equation}
Hence,
\begin{equation}
\sum_{ \substack{j\in J\setminus J_{K}} }  \Phi(B_{t_{j}},t_{j}) (B_{t_{j+1}}-B_{t_{j}})\stackrel{P}{\longrightarrow } \int_0^T \Phi_{t} d B_{t} , \ \ as \|\bigtriangleup_{n}\| \longrightarrow 0.
\end{equation}

(2) For $ \Phi_{t} \in \mathscr{L}_{2}$, $\exists \Phi_{t}^{(n)} \in \mathscr{L}_{0}$, such that $||\Phi_{t}^{(n)}-\Phi||_{2} \longrightarrow 0$, $n\longrightarrow \infty$ (see [13] Proposition1.5 in pp18 , pp26-27), according to above conclusion,

\begin{equation}
\begin{array}{ll}
&\displaystyle \int_{0}^{T}\Phi_{t}d B_{t} \equiv (\mu_{2}) \lim_{n\longrightarrow \infty} \int_{0}^{T} \Phi_{t}^{(n)}d B_{t} \\
 &=\displaystyle \mathop{l.i.m}\limits_{\|\bigtriangleup_{n}\| \longrightarrow 0} \sum_{ \substack{j\in J}}  \Phi^{(n)}(B_{t_{j}},t_{j}) (B_{t_{j+1}}-B_{t_{j}}) \\
 &=\displaystyle \mathop{l.i.m}\limits_{\|\bigtriangleup_{n}\| \longrightarrow 0} \sum_{ \substack{j\in J\setminus J_{K}} }  \Phi^{(n)}(B_{t_{j}},t_{j}) (B_{t_{j+1}}-B_{t_{j}})
\end{array}
\end{equation}

(3) For $ \Phi_{t} \in \mathscr{L}_{2}^{loc}$, $\exists$ stopping time $\tau_{n}$ of $\mathscr{F}_{t}$, such that
  $\Phi_{t}^{(n)} \equiv \Phi I_{[0,\tau_{n}]} (t) \in \mathscr{L}_{2}$ (When $m<n$,  $\Phi_{t}^{(m)} = \Phi^{n} I_{[0,\tau_{m}]} (t)$.), and $||\Phi_{t}^{(n)}-\Phi||_{2}^{loc} \longrightarrow 0$, $n\longrightarrow \infty$ (see [13]Lemma 1.5,1.6 in pp14-17, pp 33).   As shown in above conclusion,

\begin{equation}
\begin{array}{ll}
&\displaystyle \int_{0}^{T}\Phi_{t}d B_{t}= \sum_{k=1}^{\infty} \int_{\tau_{k-1}\wedge T}^{\tau_{k}\wedge T} \Phi^{(k)} dB_{t}, \ \ (\tau_{0}=0) \\
&=\displaystyle \sum_{k=1}^{\infty} \sum_{ \substack{j\in [\tau_{k-1}\wedge T,\tau_{k}\wedge T]\bigcap \{ J\setminus J_{K}\} } } \Phi^{(k)} (B_{t_{j+1}}-B_{t_{j}})
\end{array}
\end{equation}
where $J$ is the index number of partition [0,T], and $J_{K}$ means abandoning finite $K$ items in $J$.  According to the above conclusion, the Ito stochastic integral still holds in the incomplete Riemann-Stieltjes sum in mean square convergence .
Which end of the proof.

\hfill $\Box$

\subsection{Incomplete Riemann-Stieltjes sum approximating of mean square integral}

In this section, we discuss the incomplete Riemann-Stieltjes sum approximating theorem of mean square integral.

\textbf{Theorem 2.} Assume $X(\cdot,t)\in L_2=\{X(\omega,t): \forall \omega\in \Omega,\forall t\in [0,T], E[X(\cdot,t)]^2<+\infty \}$ satisfying $ E[X(\cdot,t)]^2 $ is bounded or Riemann integrable, then

\begin{equation}
\begin{array}{ll}
\displaystyle \int_{0}^{T} X(t) dt&=\displaystyle \mathop{l.i.m}\limits_{\|\bigtriangleup_{n}\| \longrightarrow 0 } \sum_{j=0}^{n-1} X(u_{j}) (t_{j+1}-t_{j}) \\
                    &=\displaystyle \mathop{l.i.m}\limits_{\|\bigtriangleup_{n}\| \longrightarrow 0 } \sum_{\substack{j\in J\setminus J_{K}}} X(u_{j}) (t_{j+1}-t_{j}) \\
\end{array}
\end{equation}

\textbf{Proof}  $L_2=\{X(\cdot,t): E[X(\cdot,t)]^2<+\infty, \forall t\in [0,T]\}$ is a linear space. Definite the inner product on $L_2$ as
\begin{equation}
(X,Y)=E[XY],\ \ \forall X,Y \in L_2.
\end{equation}

For $\forall X \in L_2$, denote $\| X \|=(X,X)^{\frac{1}{2}}$. That is $\| X \|=(E[X^2])^{\frac{1}{2}}$. Then, $L_2$ is a complete inner product space (Hilbert space)[8,17]. According to distance--inequality, we have
\begin{equation}
\begin{array}{ll}
\displaystyle E|\sum_{\substack{j\in J_{K}}} X(u_{j}) (t_{j+1}-t_{j})|^2 \\
=\displaystyle \| \sum_{\substack{j\in J_{K}}} X(u_{j}) (t_{j+1}-t_{j}) \|^2 \\
\leq \displaystyle [\sum_{\substack{j\in J_{K}}} \|  X(u_{j})\| (t_{j+1}-t_{j}) ]^2 \\
\end{array}
\end{equation}

If $E[X(t)]^2 $ is bounded on [0,T],then $\exists M>0$ such that $E[X(t)]^2 <M$, as $\|  X(u_{j})\|=E^{\frac{1}{2}}[X(u_{j})]^2$, then
\begin{equation*}
\begin{array}{ll}
&\displaystyle \sum_{\substack{j\in J_{K}}} \|X(u_{j})\|(t_{j+1}-t_{j}) \\
&\leq\displaystyle MK\|\bigtriangleup_{n}\|\longrightarrow 0,\ \ \mbox{as} \ \ {\|\bigtriangleup_{n}\|\longrightarrow 0} \\
\end{array}
\end{equation*}
Then,
\begin{equation*}
\displaystyle \mathop{l.i.m}\limits_{\|\bigtriangleup_{n}\| \longrightarrow 0  } \sum_{\substack{j\in J_{K}}} X(u_{j}) (t_{j+1}-t_{j})=0
\end{equation*}

If $E[X(t)]^2$ is Riemann integrable, then $\int_{0}^{T} E|X(t)|^2 dt $ exists, and  $\int_{0}^{T} \|X(t)\|dt=\int_{0}^{T} E^{\frac{1}{2}}|X(t)|^2 dt $ is Riemann integrable, according to Theorem 2 in [14], then

\begin{equation}
\displaystyle \lim_{\|\bigtriangleup_{n}\| \longrightarrow 0  } \sum_{\substack{j\in J_{K}}} \|  X(u_{j})\| (t_{j+1}-t_{j})=0
\end{equation}

Therefore,
\begin{equation}
\lim_{\|\bigtriangleup_{n}\| \longrightarrow 0  } E|\sum_{\substack{j\in J_{K}}} X(u_{j}) (t_{j+1}-t_{j})|^2 =0.
\end{equation}

That means
\begin{equation}
\mathop{l.i.m}\limits_{\|\bigtriangleup_{n}\| \longrightarrow 0 } \sum_{\substack{j\in J_{K}}} X(u_{j}) (t_{j+1}-t_{j}) =0.
\end{equation}

According to Lemma 1, we have
\begin{equation}
\mathop{l.i.m}\limits_{\|\bigtriangleup_{n}\| \longrightarrow 0 } \sum_{j=0}^{n-1} X(u_{j}) (t_{j+1}-t_{j})
=\mathop{l.i.m}\limits_{\|\bigtriangleup_{n}\| \longrightarrow 0 } \sum_{\substack{j\in J\setminus J_{K}}} X(u_{j}) (t_{j+1}-t_{j}).
\end{equation}

\hfill $\Box$

\subsection{Incomplete Riemann-Stieltjes sum approximating of Stratonovich Stochastic integral}

We will discuss the incomplete Riemann-Stieltjes sum of Stratonovich Stochastic integral in this section, and we only investigate the form of formula (3).

Although $\Phi_{t}$ in Ito stochastic integral is nonanticipative, $\Phi(\frac{B_{t_{j}}+B_{t_{j+1}}}{2},t_{j})$ is not non-anticipating with respect to $B_{t_{j+1}}-B_{t_{j}}$. Generally $E[\Phi(\frac{B_{t_{j}}+B_{t_{j+1}}}{2},t_{j})$ $(B_{t_{j+1}}-B_{t_{j}})]\neq 0$. However, we still have following conclusion. We still follow the original idea that Stratonovich Stochastic integral is obtained by substituting $\Phi(B_{t_{j}},t_{j})$ with $\Phi(\frac{B_{t_{j}}+B_{t_{j+1}}}{2},t_{j})$ in Ito Stochastic integral.

\textbf{Lemma 6.} Suppose that $\Phi_{t}\in \mathscr{L}_{0}$, or $\Phi_{t}\in \mathscr{L}_{2}$ with $E[\Phi_{t}]^4$ is bounded on [0,T], then

\begin{equation}
\displaystyle \lim\limits_{\|\bigtriangleup_{n}\| \longrightarrow 0 } \mathbb{E} | \sum_{ \substack{j\in J_{K}} }  \Phi(\frac{B_{t_{j}}+B_{t_{j+1}}}{2},t_{j}) (B_{t_{j+1}}-B_{t_{j}}) |^2 =0
\end{equation}

\textbf{Proof } (1) If $\Phi_{t}\in \mathscr{L}_{0}$, then
\begin{equation*}
\begin{array}{ll}
&\displaystyle \mathbb{E} ^{\frac{1}{2}} [| \sum_{ \substack{j\in J_{K}} }  \Phi(\frac{B_{t_{j}}+B_{t_{j+1}}}{2},t_{j}) (B_{t_{j+1}}-B_{t_{j}}) |^2 ]\\
\leq & \displaystyle  \sum_{ \substack{j\in J_{K}} } \large{ \mathbb{E}^{\frac{1}{2}} [|\Phi(\frac{B_{t_{j}}+B_{t_{j+1}}}{2},t_{j}) (B_{t_{j+1}}-B_{t_{j}}) |^2 ] \large} \  \ (Minkowski\ \ inequality)  \\
\leq & \displaystyle  \sum_{ \substack{j\in J_{K}} } \mathbb{E}^{\frac{1}{2}} [ \sup_{\substack{j\in J_{K}}} \{\Phi^2(\frac{B_{t_{j}}+B_{t_{j+1}}}{2},t_{j})\} (B_{t_{j+1}}-B_{t_{j}})^2 ] \\
\leq & \displaystyle  \sum_{ \substack{j\in J_{K}} } \{ \sup_{\substack{j\in J_{K}}} \{\Phi^2(\frac{B_{t_{j}}+B_{t_{j+1}}}{2},t_{j})\} \}^{\frac{1}{2}}  \mathbb{E}^{\frac{1}{2}} [(B_{t_{j+1}}-B_{t_{j}})^2 ] \\
\end{array}
\end{equation*}
\begin{equation}
\begin{array}{ll}
=&\displaystyle  \sum_{ \substack{j\in J_{K}} }  \sup_{\substack{j\in J_{K}}} \{\Phi^2(\frac{B_{t_{j}}+B_{t_{j+1}}}{2},t_{j})\}^{\frac{1}{2}} {\vartriangle t_{j}}^{\frac{1}{2}}\\
\leq &\displaystyle  \sup_{\substack{j\in J_{K}}} \{\Phi^2(\frac{B_{t_{j}}+B_{t_{j+1}}}{2},t_{j})\}^{\frac{1}{2}} \sum_{ \substack{j\in J_{K}}} \{\vartriangle t_{j}\}^{\frac{1}{2}} \\
\leq &\displaystyle  \sup_{\substack{j\in J_{K}}} \{\Phi^2(\frac{B_{t_{j}}+B_{t_{j+1}}}{2},t_{j})\} K  \{\max{ \{\vartriangle t_{j}\}\}^{\frac{1}{2}}} \longrightarrow 0, \ \ \mbox{as} \ \ {\max{ \{\vartriangle t_{i}\}}\longrightarrow 0} \\
\end{array}
\end{equation}
Hence,
\begin{equation}
\displaystyle \lim\limits_{\|\bigtriangleup_{n}\| \longrightarrow 0 } \mathbb{E} [| \sum_{ \substack{j\in J_{K}} }  \Phi(\frac{B_{t_{j}}+B_{t_{j+1}}}{2},t_{j}) (B_{t_{j+1}}-B_{t_{j}}) |^2 ]=0
\end{equation}
which ends the proof.

(2) If $\Phi_{t}\in \mathscr{L}_{2}$ with $E[\Phi_{t}]^4$ is bounded,  $\exists M>0 $ so that $E[\Phi_{t}]^4 < M$, and $E[B_{t_{j+1}}-B_{t_{j}}]^4=3(t_{j+1}-t_{j})^2$, then
\begin{equation}
\begin{array}{ll}
&\displaystyle \mathbb{E} ^{\frac{1}{2}} [| \sum_{ \substack{j\in J_{K}} }  \Phi(\frac{B_{t_{j}}+B_{t_{j+1}}}{2},t_{j}) (B_{t_{j+1}}-B_{t_{j}}) |^2 ]\\
\leq & \displaystyle  \sum_{ \substack{j\in J_{K}} } { \mathbb{E}^{\frac{1}{2}} [|\Phi(\frac{B_{t_{j}}+B_{t_{j+1}}}{2},t_{j}) (B_{t_{j+1}}-B_{t_{j}}) |^2 ] } \  \ (\mbox{Minkowski\ \ inequality})  \\
 =& \displaystyle  \sum_{ \substack{j\in J_{K}} } \{\mathbb{E} |\Phi^2(\frac{B_{t_{j}}+B_{t_{j+1}}}{2},t_{j}) (B_{t_{j+1}}-B_{t_{j}})^2| \}^{\frac{1}{2}} \\
\leq & \displaystyle  \sum_{ \substack{j\in J_{K}} } \{\mathbb{E}^{\frac{1}{2}} [\Phi^4(\frac{B_{t_{j}}+B_{t_{j+1}}}{2},t_{j})] \mathbb{E}^{\frac{1}{2}}[(B_{t_{j+1}}-B_{t_{j}})^4] \}^{\frac{1}{2}} \  \ (\mbox{H$\ddot{o}$lder \ \ inequality} )\\
\leq & \displaystyle  \sum_{ \substack{j\in J_{K}} } \mathbb{E}^{\frac{1}{4}} |\Phi^4(\frac{B_{t_{j}}+B_{t_{j+1}}}{2},t_{j})| \mathbb{E}^{\frac{1}{4}}[(B_{t_{j+1}}-B_{t_{j}})^4] \\
\leq & \displaystyle  M^{\frac{1}{4}} 3^{\frac{1}{4}} \sum_{ \substack{j\in J_{K}} } {\vartriangle t_{j}}^{\frac{1}{2}}\\
\leq & \displaystyle  M^{\frac{1}{4}} 3^{\frac{1}{4}} K  \{\max{ \{\vartriangle t_{j}\}\}^{\frac{1}{2}}} \longrightarrow 0, \ \ \mbox{as} \ \ {\max{ \{\vartriangle t_{i}\}}\longrightarrow 0} \\
\end{array}
\end{equation}
Hence,
\begin{equation}
\displaystyle \lim\limits_{\|\bigtriangleup_{n}\| \longrightarrow 0 } \mathbb{E} [| \sum_{ \substack{j\in J_{K}} } \Phi(\frac{B_{t_{j}}+B_{t_{j+1}}}{2},t_{j}) (B_{t_{j+1}}-B_{t_{j}}) |^2 ]=0
\end{equation}

\hfill $\Box$

\textbf{Theorem 3.} Suppose that $\Phi_{t}\in \mathscr{L}_{0}$, or $\Phi_{t}\in \mathscr{L}_{2}$ with $E[\Phi_{t}]^4$ is bounded on [0,T], then
\begin{equation}
\displaystyle \mathop{l.i.m}\limits_{\|\bigtriangleup_{n}\| \longrightarrow 0 }  \sum_{ \substack{j\in J\setminus J_{K}} }  \Phi(\frac{B_{t_{j}}+B_{t_{j+1}}}{2},t_{j}) (B_{t_{j+1}}-B_{t_{j}})= \int_0^T \Phi_{t}\circ d B_{t}.
\end{equation}

\textbf{Proof} According to Lemma 6, we obtain
\begin{equation}
\begin{array}{ll}
&\displaystyle \int_0^T \Phi_{t}\circ d B_{t}=\displaystyle \mathop{l.i.m}\limits_{\|\bigtriangleup_{n}\| \longrightarrow 0 }\sum_{ \substack{j\in J} }  \Phi(\frac{B_{t_{j}}+B_{t_{j+1}}}{2},t_{j}) (B_{t_{j+1}}-B_{t_{j}})\\
&= \displaystyle \mathop{l.i.m}\limits_{\|\bigtriangleup_{n}\| \longrightarrow 0 }  \sum_{ \substack{j\in J\setminus J_{K}} }  \Phi(\frac{B_{t_{j}}+B_{t_{j+1}}}{2},t_{j}) (B_{t_{j+1}}-B_{t_{j}}) \\
 &\ \ +\displaystyle \mathop{l.i.m}\limits_{\|\bigtriangleup_{n}\| \longrightarrow 0 }  \sum_{ \substack{j\in J_{K}} }  \Phi(\frac{B_{t_{j}}+B_{t_{j+1}}}{2},t_{j}) (B_{t_{j+1}}-B_{t_{j}}).\\
&= \displaystyle \mathop{l.i.m}\limits_{\|\bigtriangleup_{n}\| \longrightarrow 0 }  \sum_{ \substack{j\in J\setminus J_{K}} }  \Phi(\frac{B_{t_{j}}+B_{t_{j+1}}}{2},t_{j}) (B_{t_{j+1}}-B_{t_{j}}) \\
\end{array}
\end{equation}

\hfill $\Box$

\textbf{Remark} For the Stratonovich Stochastic integral in formula (4), we can easily prove Lemma 6 and Theorem 3 still hold.


Considering Theorem 1 ,Theorem 2 and Theorem 3, we can reach the conclusion about incomplete Riemann-Stieltjes sums approaching theorems of the three stochastic integrals as follows.

\textbf{Theorem 4.} Suppose $\Phi(x,t)$ is continuous in $t$, and has the continuous derivative $\displaystyle \frac{\partial \Phi(x,t)}{\partial x} $. $\Phi(\cdot,t)$ is a non-anticipating stochastic process. $\Phi_{t}\in \mathscr{L}_{0}$, or $\Phi_{t}\in \mathscr{L}_{2}$ with $E[\Phi_{t}]^4$ is bounded on [0,T].$\displaystyle \frac{\partial \Phi(B(t),t)}{\partial B(t)} \in \mathscr{L}_{2}$ and $\displaystyle E [\frac{\partial \Phi(B(t),t)}{\partial B(t)}]^2 $ is bounded or Riemann integral. If all $\displaystyle \int_{0}^{T} \Phi(B(t),t)\circ dB(t)$ , \\
$\displaystyle \int_{0}^{T} \Phi(B(t),t) dB(t)$ and $\displaystyle \int_{0}^{T}\frac{\partial \Phi(B(t),t)}{\partial B(t)}dt$ exist, then
\begin{equation}
\begin{array}{ll}
\displaystyle \mathop{l.i.m}\limits_{\|\bigtriangleup_{n}\| \longrightarrow 0 } \sum_{\substack{j\in J\setminus J_{K}}} \Phi(\frac{B_{t_{j}}+B_{t_{j+1}}}{2},t_{j})(B_{t_{j+1}}-B_{t_{j}})\\
=\displaystyle \mathop{l.i.m}\limits_{\|\bigtriangleup_{n}\| \longrightarrow 0 } \sum_{\substack{j\in J\setminus J_{K}}} \Phi(B_{t_{j}},t_{j})(B_{t_{j+1}}-B_{t_{j}})
+\displaystyle \mathop{l.i.m}\limits_{\|\bigtriangleup_{n}\| \longrightarrow 0 } \sum_{\substack{j\in J\setminus J_{K}}} \frac{\partial \Phi(B(u_{j}),u_j)}{\partial B(t)}(t_{j+1}-t_j)
\end{array}
\end{equation}
where $\forall u_{j} \in [t_j,t_{j+1}]$.

\textbf{Proof} According to Theorem 1, Theorem 2 and Theorem 3, it is easily to prove that the conclusion holds.

\hfill $\Box$

\textbf{Remark.} In discussion of Ito stochastic integral, we do not strictly discriminate $\mathscr{L}_{p}$ from $\mathscr{L}_{p,T}$ in order to coincide with the expression in references, and so on.

\section{Extension of Incomplete Riemann-Stieltjes sum approximating Theorems}

In Section 3, for fixed $K$, the incomplete Riemann-Stieltjes sum approximating Theorems are investigated with fixed $J_{K}$ from $J$, if all possible combination of $J_{K}\subset J$ are considered, there will be a large number of incomplete Riemann-Stieltjes sums converge to correspondent stochastic integrals respectively.

Furthermore, let $K$ changes along with $n$, denote as $K(n) \in N^{+} (0<K(n)<n)$, and $J_{K(n)} =\{i_1 , \cdots ,i_{K(n)} \} \subset J $. And, $J\setminus J_{K(n)}=\{0,1,\cdots,n-1 \}\setminus \{i_1 , \cdots ,i_{K(n)} \}$.
Since the partition $\bigtriangleup_{n}$ of [0,T] could be arbitrary, to avoid complicated discussion, we only discuss the concise case of equal partition. which means $\Delta t_{j}=\displaystyle \frac{T}{n}, j=0,1,\cdots,n-1$, and $\|\bigtriangleup_{n}\|=\displaystyle \frac{T}{n}$. We can extend all of Theorem 1,Theorem 2,Theorem 3 and Theorem 4.

\textbf{Theorem 5. } If $\int_0^T \Phi_{t} d B_{t}$ exists, $\bigtriangleup_{n}$ is an equal partition of [0,T],and $\displaystyle \lim\limits_{\substack{n\longrightarrow +\infty}} \frac{K(n)}{n}=0$, then
\begin{equation}
\displaystyle \mathop{l.i.m}\limits_{\|\bigtriangleup_{n}\| \longrightarrow 0 } \sum_{ \substack{j\in J\setminus J_{K(n)}} }  \Phi(B_{t_{j}},t_{j}) (B_{t_{j+1}}-B_{t_{j}})= \int_0^T \Phi_{t} d B_{t}.
\end{equation}

\textbf{Proof}: We only examine the formula (24), and it is easily to see formula(24) holds for $\mathscr{L}_{0}$ when $K$ is replaced by $K(n)$, then it does hold for $\mathscr{L}_{2}$ and $\mathscr{L}_{2}^{loc}$ subsequently, which ends the proof.

\hfill $\Box$

\textbf{Theorem 6.} If $\int_{0}^{T} X(t)dt$ exists, $\bigtriangleup_{n}$ is an equal partition of [0,T],and $\displaystyle \lim\limits_{\substack{n\longrightarrow +\infty}} \frac{K(n)}{n}=0$, then
\begin{equation}
\displaystyle \mathop{l.i.m}\limits_{\|\bigtriangleup_{n}\| \longrightarrow 0 } \sum_{\substack{j\in J\setminus J_{K(n)}}} X(u_{j})(t_{j+1}-t_{j})=\int_{0}^{T} X(t)dt
\end{equation}
where $\forall u_{j}\in [t_{j},t_{j+1}]$.

\textbf{Proof}: The proof is similar to Theorem 5, it only need to check formula (35)(36) and Theorem 5 in [14] when $K$ is replaced by $K(n)$.

\hfill $\Box$

\textbf{Theorem 7.} If $\int_{0}^{T} \Phi(B_{t},t) \circ dB_{t}$ exists, $\bigtriangleup_{n}$ is an equal partition of [0,T],and $\displaystyle \lim\limits_{\substack{n\longrightarrow +\infty}} \frac{K(n)}{\sqrt{n}}=0$, then
\begin{equation}
\mathop{l.i.m}\limits_{\|\bigtriangleup_{n}\| \longrightarrow 0 } \sum_{\substack{j\in J\setminus J_{K(n)}}} \Phi(\frac{B_{t_{j}}+B_{t_{j+1}}}{2},t_{j})[B_{t_{j+1}}-B_{t_{j}}]=\int_{0}^{T} \Phi(B_{t},t) \circ dB_{t}
\end{equation}

\textbf{Proof}: The proof is similar to Theorem 5, we can only check  formula (41) and (43) when $K$ is replaced by $K(n)$. They still hold. Hence it ends the proof.

\hfill $\Box$

\textbf{Theorem 8.} Suppose all the conditions of Theorem 4 are satisfied. $\bigtriangleup_{n}$ is an equal partition of [0,T], and $\displaystyle \lim\limits_{\substack{n\rightarrow +\infty}} \frac{K(n)}{\sqrt{n}}=0$, then
\begin{equation}
\begin{array}{ll}
\displaystyle \mathop{l.i.m}\limits_{\|\bigtriangleup_{n}\|\rightarrow 0}\sum_{\substack{j\in J\setminus J_{K(n)}}} \Phi(\frac{B_{t_{j}}+B_{t_{j+1}}}{2},t_{j})(B_{t_{j+1}}-B_{t_{j}})\\
=\displaystyle \mathop{l.i.m}\limits_{\|\bigtriangleup_{n}\|\rightarrow 0}\sum_{\substack{j\in J\setminus J_{K(n)}}} \Phi(B_{t_{j}},t_{j})(B_{t_{j+1}}-B_{t_{j}})
+\displaystyle \mathop{l.i.m}\limits_{\|\bigtriangleup_{n}\|\rightarrow 0}\sum_{\substack{j\in J\setminus J_{K(n)}}} \frac{\partial \Phi(B(u_{j}),u_j)}{\partial B(t)}(t_{j+1}-t_j)
\end{array}
\end{equation}
where $\forall u_{j} \in [t_j,t_{j+1}]$.


For fixed natural number $K\in N^{+} (0<K<n)$, $J_K=\{i_1,\cdots,i_K\} \subset J$ is one possible combination, there will be $C_{n}^{K}$ incomplete Riemann-Stieltjes sums. Therefore, we have the following theorems.

\textbf{Theorem 9. } For fixed K,
\begin{equation}
\begin{array}{l}
\displaystyle \# \{\sum_{\substack{j\in J\setminus J_{k}}}\Phi(B_{t_{j}},t_{j})(B_{t_{j+1}}-B_{t_{j}}) | \forall J_{k}=\{i_1,\cdots,i_k\}\subset J, 1\leq k\leq K\} \\
=C_{n}^{K}(2^{K}-1).
\end{array}
\end{equation}
For fixed K(n),
\begin{equation}
\begin{array}{l}
\displaystyle \# \{\sum_{\substack{j\in J\setminus J_{k}}}\Phi(B_{t_{j}},t_{j})(B_{t_{j+1}}-B_{t_{j}}) | \forall J_{k}=\{i_1,\cdots,i_k\}\subset J, 1\leq k\leq K(n)\} \\
=C_{n}^{K(n)} (2^{K(n)}-1).
\end{array}
\end{equation}

\textbf{Theorem 10. } For fixed K,
\begin{equation}
\begin{array}{l}
\displaystyle \# \{\sum_{\substack{j\in J\setminus J_{k}}} X(u_{j})(t_{j+1}-t_{j})|\forall J_{k}=\{i_1,\cdots,i_k\}\subset J, 1\leq k\leq K\} \\
=C_{n}^{K}(2^{K}-1).
\end{array}
\end{equation}
For fixed K(n),
\begin{equation}
\begin{array}{l}
\displaystyle \# \{ \sum_{ \substack{j\in J\setminus J_{k}}} X(u_{j})(t_{j+1}-t_{j})|\forall J_{k}=\{i_1,\cdots,i_k\}\subset J, 1\leq k\leq K(n) \} \\
=C_{n}^{K(n)} (2^{K(n)}-1).
\end{array}
\end{equation}
where $\forall u_{j}\in [t_{j},t_{j+1}] $.

\textbf{Theorem 11. } For fixed K,
\begin{equation}
\begin{array}{l}
\displaystyle \# \{\sum_{\substack{j\in J\setminus J_{k}}} \Phi(\frac{B_{t_{j}}+B_{t_{j+1}}}{2},t_{j})[B_{t_{j+1}}-B_{t_{j}}]|\forall J_{k}=\{i_1,\cdots,i_k\}\subset J,1\leq k\leq K \} \\
=C_{n}^{K}(2^{K}-1).
\end{array}
\end{equation}
For fixed K(n),
\begin{equation}
\begin{array}{l}
\displaystyle \# \{\sum_{\substack{j\in J\setminus J_{k}}} \Phi(\frac{B_{t_{j}}+B_{t_{j+1}}}{2},t_{j})[B_{t_{j+1}}-B_{t_{j}}]| \forall J_{k}=\{i_1,\cdots,i_k\}\subset J,1\leq k\leq K(n)\}\\
=C_{n}^{K(n)}(2^{K(n)}-1).
\end{array}
\end{equation}

\textbf{Remark.} Specifically,when $K(n)=n^r$, $0<r<1$, $\displaystyle \lim_{\substack{n\longrightarrow +\infty}} \displaystyle \frac{K(n)}{n}=0$,and $\displaystyle \lim_{\substack{n\longrightarrow +\infty}} 2^{K(n)}-1=\displaystyle \lim_{\substack{n\longrightarrow +\infty}} 2^{n^r}-1= +\infty$. According to Theorem 1,Theorem 2,Theorem 5 and Theorem 6, all the incomplete Riemann-Stieltjes sums in Theorem 9 and Theorem 10 approach Ito  stochastic integral and mean square integral respectively. This means there will be a large number of incomplete  Riemann-Stieltjes sums which will converge to the Ito  stochastic integral and mean square integral respectively in mean square convergence sense.

For Stratonovich stochastic integral , if $K(n)=n^r$, $0<r<1/2$, then $\displaystyle \lim_{\substack{n\longrightarrow +\infty}} \displaystyle \frac{K(n)}{n}=0$, and  $2^{K(n)}-1=2^{n^r}-1\rightarrow +\infty$, as $n\rightarrow +\infty$ still hold. According to Theorem 3 and Theorem 7, all the incomplete Riemann-Stieltjes sums in Theorem 11 will converge to Stratonovich stochastic integral.

In addition, the $J_K$s in Theorem 4 can not be same, so do the $J_{K(n)}$s in Theorem 8.

Furthermore, when $r\in (0,1)$ for incomplete Riemann-Stieltjes sums of Ito  stochastic integral and mean square integral, and
$r\in (0,\displaystyle \frac{1}{2})$ for Stratonovich stochastic integral are considered, the number of incomplete Riemann-Stieltjes sums for each of the three stochastic integrals is uncountable, while $n\rightarrow +\infty$.

Even in equal partition case, we construct a lot of Riemann-Stieltjes sums to approach each stochastic integral in mean square convergence sense. For the arbitrary partition cases, there ought to exist more numerous incomplete Riemann-Stieltjes sums converging to the specified stochastic integral. Hence, there exist uncountable infinite number of  Riemann-Stieltjes sums while the complete Riemann-Stieltjes sums converge to each stochastic integral.

According to above analyses, we obtain the following theorem.

\textbf{Theorem 12. } There are uncountable incomplete Riemann-Stieltjes sums for Ito  stochastic integral, mean square integral and Stratonovich stochastic integral respectively converging to their own stochastic integral in mean square convergence sense, along with $n\rightarrow +\infty$.

\textbf{Proof.} Since Ito stochastic integral, mean square integral and Stratonovich stochastic integral are integrable, we only check the equal partition case. According to Theorem 9,Theorem 10 and Theorem 11, the number of incomplete Riemann--Stieltjes sums with respect to $K(n)$ can be calculated, given fixed n. And, even for $r\in (0,1)$ or $r\in (0,\displaystyle \frac{1}{2})$ \begin{equation}
\displaystyle \lim_{\substack{n\longrightarrow +\infty}} C_{n}^{K(n)} (2^{K(n)}-1)=\displaystyle \lim_{\substack{n\longrightarrow +\infty}} C_{n}^{n^r} (2^{n^r}-1)= +\infty
\end{equation}
To investigate whether the three  incomplete Riemann--Stieltjes sum sets with respect to $K(n)$  in Theorem 9,Theorem 10 and Theorem 11 are countable or uncountable while $n\rightarrow +\infty$, we now consider the following sets along with the process that $n\rightarrow +\infty$ as $\|\bigtriangleup_{n}\| \rightarrow 0$.
\begin{equation}
\begin{array}{l}
\displaystyle \# \{\sum_{\substack{j\in J\setminus J_{k}}}\Phi(B_{t_{j}},t_{j})(B_{t_{j+1}}-B_{t_{j}})|\forall J_{k}=\{i_1,\cdots,i_k\}\subset J, 1\leq k\leq K(n)=n^r, \\
 \  \   \    \     \  \ \forall r\in(0,1), \displaystyle \lim_{\substack{n\rightarrow +\infty}} \displaystyle \frac{K(n)}{n}=0 \} \\
\end{array}
\end{equation}
\begin{equation}
\begin{array}{l}
\displaystyle \# \{\sum_{\substack{j\in J\setminus J_{k}}}X(u_{j})(t_{j+1}-t_{j})|\forall J_{k}=\{i_1,\cdots,i_k\}\subset J, 1\leq k\leq K(n)=n^r, \\
 \  \   \    \     \  \ \forall r\in(0,1), \displaystyle \lim_{\substack{n\rightarrow +\infty}} \displaystyle \frac{K(n)}{n}=0 \} \\
\end{array}
\end{equation}
where $\forall u_{j}\in [t_{j},t_{j+1}] $.
\begin{equation}
\begin{array}{l}
\displaystyle \# \{\sum_{\substack{j\in J\setminus J_{k}}} \Phi(\frac{B_{t_{j}}+B_{t_{j+1}}}{2},t_{j})[B_{t_{j+1}}-B_{t_{j}}]| \forall J_{k}=\{i_1,\cdots,i_k\}\subset J,1\leq k\leq K(n) \\
 \  \   \    \     \  \ =n^r,\forall r\in(0,\frac{1}{2}), \displaystyle \lim_{\substack{n\rightarrow +\infty}} \displaystyle \frac{K(n)}{\sqrt{n}}=0 \} \\
\end{array}
\end{equation}
We only discuss the case $r\in(0,\frac{1}{2})$. The above three sets are related to natural number $n$ and real number $r$. Since $n^{r}$ is strictly monotonic increasing function in $r\in(0,\frac{1}{2})$. For $\forall 0<r_1<r_2<\frac{1}{2}$ , whatever how much $r_1, r_2$ could be small, as $n^{r_1}<n^{r_2}$, and $\displaystyle \lim_{\substack{n\rightarrow +\infty}} n^{r_1}=\displaystyle \lim_{\substack{n\rightarrow +\infty}} n^{r_2}=+\infty$, there exist $n$ which is large enough such that
$[n^{r_1}]<< [n^{r_2}]$, here $[\cdot]$ is the integer function. Under this equal partitions $\bigtriangleup_{n}$, $\{J\setminus J_{[n^{r_1}]}\}$ and $\{J\setminus J_{[n^{r_2}]}\}$ could be different, then the incomplete Riemann--Stieltjes sums expressions in each of above three sets are different with different $r_1$ and $r_2$. Hence, we construct a mapping function from $r\in(0,\frac{1}{2})$ to each subset of formula (59)(60)(61). Since $(0,\frac{1}{2})$ is uncountable, the numbers of the three sets in formula (59)(60)(61) are uncountable as $n\rightarrow +\infty$. Therefore it ends the proof.

\hfill $\Box$

\section{Simulation of Incomplete Riemann-Stieltjes sum Approximation of Stochastic Integral}

The Brownian motion involved in the simulation is generated according to algorithm in [22].

\subsection{Simulation of Incomplete Riemann-Stieltjes sum approximation of Ito stochastic integral}

To demonstrate the simulation of incomplete Riemann-Stieltjes sum approximating theorem of Ito stochastic integral
\begin{equation}
\displaystyle \int_{0}^{T} B_{t} dB_{t}=\displaystyle \frac{1}{2} B_{T}^2 -\frac{1}{2} T=\displaystyle \lim_{\substack{n\longrightarrow +\infty}} \sum_{\substack{j\in J\setminus J_{K(n)}}}  B_{t_{j}}(B_{t_{j+1}}-B_{t_{j}}).
\end{equation}

we design two simulations:

(1) Let $T=1$,$n=10^6$,$r=\{0,0.1,0.2,\cdots,0.9\}$, where $r=0$ represents the normal standard complete Riemann-Stieltjes sum, and $J=\{0,1,2,\cdots,n-1\}$.  Then, we adopt three strategies,
    $J_{n^r}=\{0,1,2,\cdots,n^r-1\}$, $J_{n^r}=\{i_1,i_2,\cdots,i_{n^r}\} \in J$, $\forall \{i_1,i_2,\cdots,i_{n^r}\} \in J$, and $J_{n^r}=\{n-n^r,\cdots,n-2,n-1\}$, which means deleting the forth $n^r$ items in complete Riemann-Stieltjes sum, randomly deleting $n^r$ items in complete Riemann-Stieltjes sum and deleting the end $n^r$ items in complete Riemann-Stieltjes sum, and denote as $IRS_b$,$IRS_r$,$IRS_e$ respectively. The simulation results are shown in Figure 1.

\begin{figure}[!htbp]
\begin{center}
\begin{tabular}{c}
\includegraphics[width=10cm]{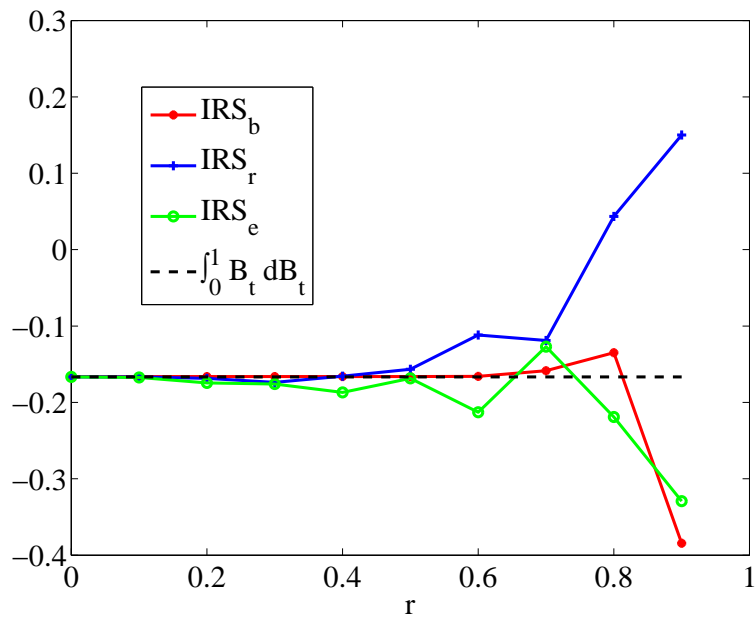} \\
\end{tabular}
\caption{Comparison of incomplete Riemann-Stieltjes sums in styles of deleting the forth $n^r$ items , randomly deleting $n^r$ items and deleting the end $n^r$ items of complete Riemann-Stieltjes sum and Ito stochastic integral}
\end{center}
\end{figure}

(2) Above experiment is only one realization of stochastic orbit. We iterate the simulation $N=10^5$ times, and define the statistic  mean absolute error (MAE)  as follows:

\begin{equation}
\displaystyle MAE=\displaystyle \frac{1}{N} \sum_{i=0}^{N} |[\sum_{\substack{j\in J\setminus J_{K(n)}}}B_{t_{j}}(B_{t_{j+1}}-B_{t_{j}})]-[\frac{1}{2} B_{T}^2 -\frac{1}{2}T]|.
\end{equation}

where $T=1$, $n=10^6$,$r=\{0,0.1,0.2,\cdots,0.9\}$, and we still adopt three strategies,
    $J_{n^r}=\{0,1,2,\cdots,n^r-1\}$, $J_{n^r}=\{i_1,i_2,\cdots,i_{n^r}\}$, $\forall \{i_1,i_2,\cdots,i_{n^r}\} \in J$, and  $J_{n^r}=\{n-n^r,n-n^r+1,\cdots,n-1\}$, , and denote as $MAE_b$,$MAE_r$,$MAE_e$ respectively. The simulation results are shown in Figure 2.

\begin{figure}[!htbp]
\begin{center}
\begin{tabular}{c}
\includegraphics[width=10cm]{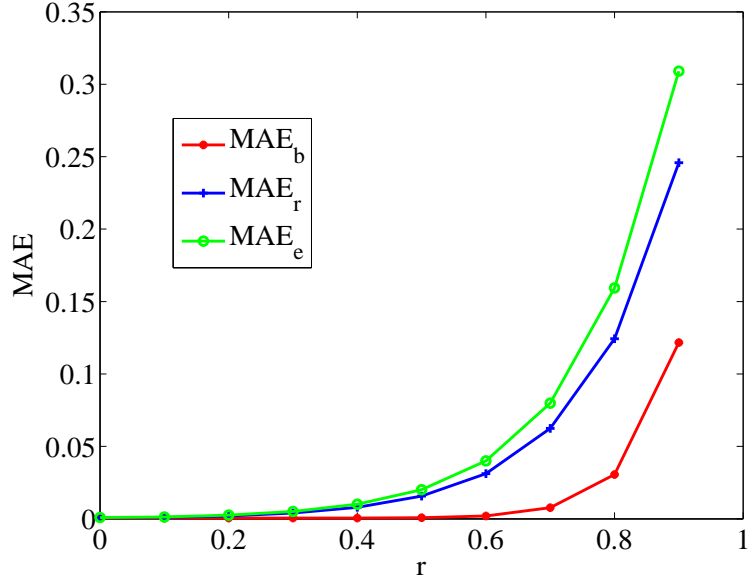} \\
\end{tabular}
\caption{MAE indexes of deleting the forth $n^r$ items, randomly deleting $n^r$ items and deleting the end $n^r$ items in complete Riemann-Stieltjes sum of Ito stochastic integral}
\end{center}
\end{figure}

Figure 1 show that the more the deleting items, the more derivation of the approximation result. And, while $r=0.1$, deleting $\displaystyle n^{\frac{1}{10}}$ item seems have no influence on simulation of $\displaystyle \int_{0}^{T} B_{t} dB_{t}$.
While increasing the simulation iterations, Figure 2 verify the above result in Figure 1, and also show that deleting the forth items strategy is better than the other two strategies. And deleting end strategy is the worse result. In fact, $\displaystyle \int_{0}^{T} B_{t} dB_{t}$ is path--dependent, deleting the end items will greatly affect the final integral result. Hence the deleting item of Ito stochastic integral is also path--dependent.

\subsection{Simulation of Incomplete Riemann-Stieltjes sum Approximation of Stratonovich Stochastic Integral}

To demonstrate the simulation of incomplete Riemann-Stieltjes sum approximating theorem of Stratonovich stochastic integral

\begin{equation}
\displaystyle \int_{0}^{T} B_{t}\circ dB_{t}=\displaystyle \frac{1}{2} B_{T}^2=\displaystyle \lim_{\substack{n\longrightarrow +\infty}} \sum_{\substack{j\in J\setminus J_{K(n)}}} \frac{B_{t_{j+1}}+ B_{t_{j}}}{2}(B_{t_{j+1}}-B_{t_{j}}).
\end{equation}

we design two simulations:

(1) Let $T=1$,$n=10^6$,$r=\{0,0.1,0.2,\cdots,0.9\}$, where $r=0$ represents the normal standard complete Riemann-Stieltjes sum, and $J=\{0,1,2,\cdots,n-1\}$. Then, we adopt three strategies,
    $J_{n^r}=\{0,1,2,\cdots,n^r-1\}$, $J_{n^r}=\{i_1,i_2,\cdots,i_{n^r}\}$, $\forall \{i_1,i_2,\cdots,i_{n^r}\} \in J$, and     $J_{n^r}=\{n-n^r,n-n^r+1,\cdots,n-1\}$, which means deleting the forth $n^r$ items in complete Riemann-Stieltjes sum, randomly deleting $n^r$ items in complete Riemann-Stieltjes sum and deleting the end $n^r$ items in complete Riemann-Stieltjes sum, and denote as $IRS_b$,$IRS_r$,$IRS_e$ respectively. The simulation results are shown in Figure 3.

\begin{figure}[!htbp]
\begin{center}
\begin{tabular}{c}
\includegraphics[width=10cm]{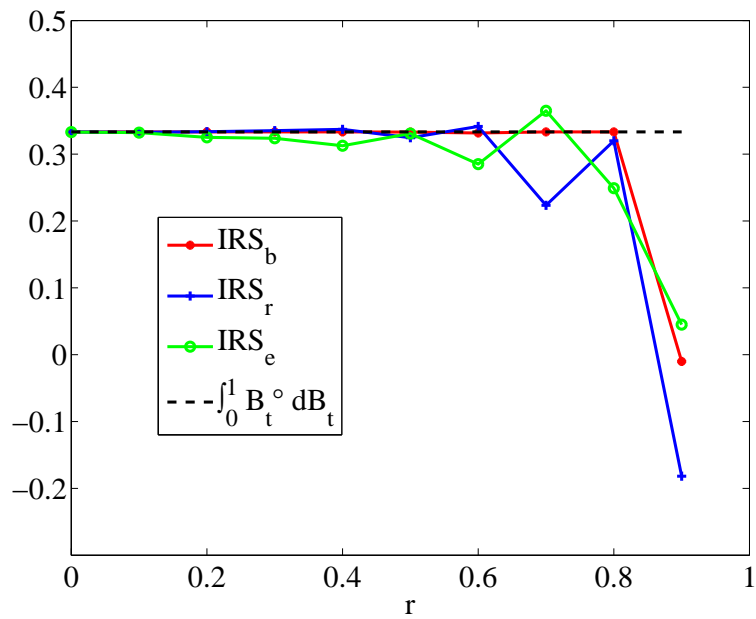}
\end{tabular}
\caption{Comparison of incomplete Riemann-Stieltjes sums in styles of deleting the forth $n^r$ items , randomly deleting $n^r$ items and deleting the end $n^r$ items of complete Riemann-Stieltjes sum and Stratonovich stochastic integral}
\end{center}
\end{figure}

(2) Above experiment is only one realization of stochastic path. We iterate the simulation $N=10^5$ times, and define the statistic mean absolute error (MAE)  as follows:

\begin{equation}
\displaystyle MAE=\displaystyle \frac{1}{N} \sum_{i=0}^{N} |[\sum_{\substack{j\in J\setminus J_{K(n)}}}\frac{B_{t_{j+1}}+ B_{t_{j}}}{2}(B_{t_{j+1}}-B_{t_{j}})]-\frac{1}{2} B_{T}^2 |.
\end{equation}

where $T=1$, $n=10^6$,$r=\{0,0.1,0.2,\cdots,0.9\}$, and we still adopt three strategies,
    $J_{n^r}=\{0,1,2,\cdots,n^r-1\}$, $J_{n^r}=\{i_1,i_2,\cdots,i_{n^r}\}$, $\forall \{i_1,i_2,\cdots,i_{n^r}\} \in J$, and $J_{n^r}=\{n-n^r,n-n^r+1,\cdots,n-1\}$, and denote as $MAE_b$,$MAE_r$,$MAE_e$ respectively. The simulation results are shown in Figure 4.

\begin{figure}[!htbp]
\begin{center}
\begin{tabular}{c}
\includegraphics[width=10cm]{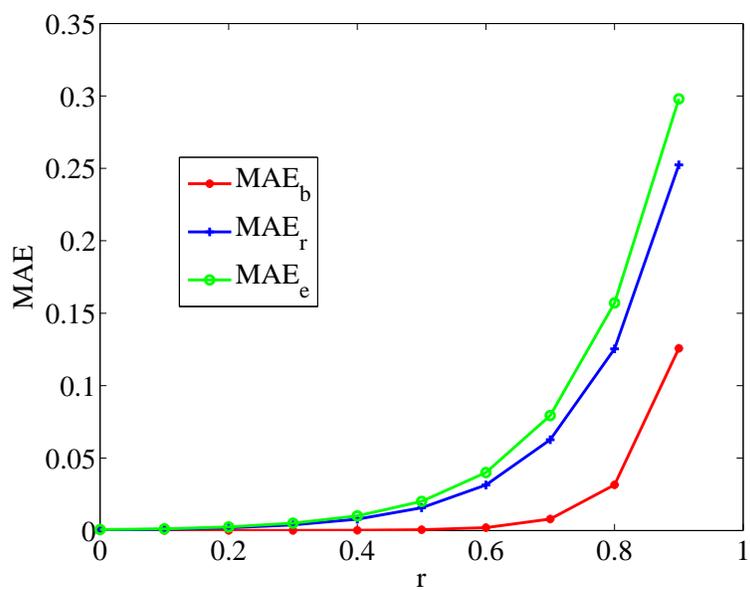} \\
\end{tabular}
\caption{MAE indexes of deleting the forth $n^r$ items, randomly deleting $n^r$ items and deleting the end $n^r$ items in complete Riemann-Stieltjes sum of Stratonovich stochastic integral}
\end{center}
\end{figure}

Note that in Theorem 7 and remark of Section 4, the sufficient condition for convergence of incomplete Riemann-Stieltjes sum of Stratonovich stochastic integral is $K(n)=n^r$, $0<r<1/2$. We perform on the $r=\{0,0.1,0.2,\cdots,0.9\}$ to show the divergence of simulation on $r\in \{0.5,0.6,\cdots,0.9\}$.
We only propose an outline of simulation, and $n=10^6$ is relatively small for $n\rightarrow +\infty$. The limit behavior simulations of incomplete Riemann-Stieltjes sums of stochastic integrals about complicated stochastic processes need high quality computer platform.

\section{Conclusion}

The incomplete Riemann-Stieltjes sum approximating theorems are developed for Ito stochastic integral, mean square integral and Stratonovich stochastic integral respectively. There are uncountable incomplete Riemann-Stieltjes sums approach to each stochastic integral when the originally defined Riemann-Stieltjes sum is converging, which reveals the complexity of simulation stochastic integral especially in stochastic finance and option pricing and non-uniqueness of stochastic integral and stochastic differential equation in stochastic process realization. In addition, incomplete Riemann-Stieltjes sum style can also be treated as a kind of discontinuity of integral interval.Hence, it provides a new view for the integral in convergence sense.  Furthermore,
our research is mainly based on Riemann framework, however Lebesgue integral is widely used in measure theory, probability theory and stochastic process, for example, the definitions of expectation,martingale and semimartingale are mostly in Lebesgue sense, how to realize and simulate the stochastic integral with respect to general martingale and semimartingale except Brownian motion, develop realizable expression for stochastic integral in Lebesgue sense and apply the theorems in stochastic finance  would be our future interests.



\begin{thebibliography}{14}

\bibitem{} K. It$\hat{o}$. Stochastic integral, Proc. Imp. Acad. Tokyo. 1944,20,pp 519-524.

\bibitem{} D.L. Fisk. Quasi-martingales and stochastic integrals. Techn. Report Dept. Math. Michigan State Univ. , 1. 1963.
\bibitem{} R.L. Stratonovich. Introduction to the theory of random noise. Gordon and Breach. New York, London. 1963.
\bibitem{} C. Gardiner. Stochastic methods: a handbook for the natural and social sciences (4th Edition). Springer-Verlag Berlin Heidelberg. 2009. pp79-98.



\bibitem{} M.H.A Davis, Martingale integrals and stochastic calculus, in Communication Systems and Random Process Theory, J.K. Skwirzynski (Ed.), Sisthoff and Nordhoff, Holland,1978.
\bibitem{} K.L. Chung, R.J. Williams. Introduction to Stochastic Integration (2nd Edition). Birkh$\ddot{a}$user.Boston,Basel,Berlin. 1990.
\bibitem{} I. Karatzas, S. E. Shreve. Brownian Motion and Stochastic Calculus (Second Edition). Springer-Verlag. 1991.
\bibitem{} P.E. Kloeden, E. Platen. Numerical solution of stochastic differential equations. Applications of Mathematics. Springer-Verlag. Berlin, New York. 1992. pp 81-100.
\bibitem{} B. $\emptyset$ ksendal. Stochastic Differential Equations: An Introduction with Application. Springer-Verlag Berlin Heidelberg New York. 2003, pp 22-37.
\bibitem{} D. Revuz, M. Yor. Continuous Martingales and Brownian Motion(3rd Edition). Springer, 2005.
\bibitem{} F.C. Klebaner. Introduction to Stochastic Calculus with Application (3rd Edition). World Scientific Publishing.2012.
\bibitem{} H. Matsumoto, S. Taniguchi. Stochastic Analysis: It$\hat{o}$ and Malliavin Calculus in Tandem. Cambridge University Press. 2017.
\bibitem{} G.L Gong. Introduction to Stochastic differential equation (2nd Edition). Peking University Press. Beijing. 1995. p 9-39. (Chinese)


\bibitem{} J.W Liu, Y. Liu. Deleting Items and Disturbing Mesh Theorems for Riemann Definite Integral and Their Applications. https://arxiv.org/abs/1702.04464.


\bibitem{} J.L. Doob. Stochastic processes. John Wiley \& Sons. New York. 1991. pp92-95.
\bibitem{} P.H. Wirsching, T.L. Paez, K. Ortiz. Random vibrations-theory and practice. John Wiley \& Sons. New York. 1995. pp409-410.
\bibitem{} Y.L. Lin. Applied Stochastic Processes. Tsinghua University Press. Beijing. 2002. pp282-288. (Chinese)
\bibitem{} M. Lo$\grave{e}$ve. Probability Theory.II(4th Edition). Springer-Verlag. New York.1991. pp130-139.

\bibitem{} T.L. Toh, T.S. Chew. The Riemann approach to stochastic integration using non-uniform meshes. J. Math. Anal. Appl. 2003. 280, pp133-147.

\bibitem{} A.A. Borovkov. Probability theory (5th Edition). Springer-Verlag London 2013. pp88.

\bibitem{} S.H. Cheng. Advanced probability theory. Peking University Press: Beijing. 1996. pp6-7. (Chinese)

\bibitem{} D.J. Higham. An Algorithmic Introduction to Numerical Simulation of Stochastic Differential Equations, SIAM Rev. (Educ. Sect.), 2001, 43(3), 525-546.


\end{thebibliography}
\end{document}